\documentclass[12pt]{article}
\begin{document}
\def\supp{\mathrm{supp}\,}
\def\C{\mathbf{C}}
\def\Z{\mathbf{Z}}
\def\Del{\nabla}
\def\P{\mathbf{P}}
\def\Rea{{\mathrm{Re}\,}}
\title{Brody curves omitting hyperplanes}
\author{Alexandre Eremenko\thanks{Supported by NSF grant
DMS-0555279.}}
\maketitle
\begin{abstract} A {\em Brody curve}, a.k.a. normal curve,
is a holomorphic map $f$ from the complex line
$\C$ to the complex projective space
$\P^n$ such that the family of its
translations $\{ z\mapsto f(z+a):a\in\C\}$ is normal.
We prove that Brody curves omitting $n$ hyperplanes
in general position have growth order at most one, normal
type. This generalizes a result of Clunie and Hayman
who proved it for $n=1$.
\end{abstract}

\noindent
{\bf Introduction}
\vspace{.1in}

We consider holomorphic curves $f:\C\to\P^n$.
The spherical derivative $\| f'\|$ measures the
length distortion from the
Euclidean metric in $\C$ to the
Fubini--Study metric in $\P^n$. The explicit expression
is 
$$\| f'\|^2=
\| f\|^{-4}\sum_{i\neq j}|f_i^\prime f_j-f_if_j^\prime|^2,$$
where $(f_0,\ldots,f_n)$ is a homogeneous
representation of $f$
(that is the $f_j$ are entire functions which
never simultaneously
vanish), and
$$\| f\|^2=\sum_{j=0}^n|f_j|^2.$$
A holomorphic curve is called a Brody curve 
if its spherical derivative is bounded.
This is equivalent to normality of the family of 
translations
$\{ z\mapsto f(z+a):a\in\C\}$.

Brody curves are important for at least two reasons.
First one is the rescaling trick known as Zalcman's lemma or
Brody's lemma: for every non-constant holomorphic curve
$f$ one can find a sequence of affine maps $a_k:\C\to\C$
such that the limit $f\circ a_k$ exists
and is a non-constant Brody curve.
Second reason is Gromov's theory of mean dimension
\cite{Gromov} in which a space of Brody curves
is one of the main examples.

For the recent work
on Brody curves we refer to \cite{E,T1,T2,T3,W}.
A general reference for holomorphic curves is \cite{Lang}.

We recall that the Nevanlinna characteristic is
defined by
$$T(r,f)=\int_0^r \frac{dt}{t}
\left(\frac{1}{\pi}\int_{|z|\leq t}\| f'\|^2(z)dm_z\right),$$
where $dm$ is the area element in $\C$.
So Brody curves have order at most two normal type, that is
\begin{equation}\label{0}
T(r,f)=O(r^2).
\end{equation}

Clunie and Hayman \cite{CH} found that Brody
curves $\C\to\P^1$ omitting one
point in $\P^1$ must have smaller order of growth:
\begin{equation}
\label{1}
T(r,f)=O(r).
\end{equation}
A different proof of this fact is due to 
Pommerenke \cite{Pom}.
In this paper we prove that this phenomenon
persists in all dimensions.
\vspace{.1in}

\noindent
{\bf Theorem.} {\em Brody curves $f:\C\to\P^n$
omitting $n$ hyperplanes in general position satisfy 
$(\ref{1})$.}
\vspace{.1in}

Under the stronger assumption that a Brody
curve omits $n+1$
hyperplanes in general position, the same conclusion
was obtained by Berteloot and Duval \cite{BD} and
Tsukamoto \cite{T2}, with different proofs.

Combined with a result of Tsukamoto \cite{T1} our theorem 
implies
\vspace{.1in}

\noindent
{\bf Corollary.} {\em Mean dimension in the sense of Gromov
of the space of Brody curves in
$$\P^n\backslash\{ n\;\mbox{hyperplanes in general position}\}$$
is zero.}
\vspace{.1in}

The condition that $n$ hyperplanes are omitted is exact:
it is easy to show by direct computation that the curve
$(f_0,f_1,1,\ldots,1)$, where $f_i$ are
appropriately chosen entire functions 
such that $f_1/f_0$ is an elliptic function,
is a Brody curve, it
omits $n-1$ hyperplanes,
and $T(r,f)\sim cr^2,\; r\to\infty$ where $c>0$. 
This example will be discussed in the end of the paper.

The author thanks Alexandr Rashkovskii and
Masaki Tsukamoto 
for inspiring
conversations on the subject.
\vspace{.2in}

\noindent
{\bf Preliminaries}
\vspace{.2in}

Without loss of generality
we assume that
the omitted hyperplanes are given in the
homogeneous coordinates
by the equations
$\{ w_j=0\},\; 1\leq j\leq n.$
We fix a homogeneous representation
$(f_0,\ldots,f_n)$ of our curve, where $f_j$
are entire functions without common zeros, and $f_n=1$.
We assume without loss of generality that $f_0(0)\neq 0$.

Then
\begin{equation}
\label{u}
u=\log\sqrt{|f_0|^2+\ldots+|f_n|^2}
\end{equation}
is a positive subharmonic function, and Jensen's formula
gives
$$T(r,f)=\frac{1}{2\pi}\int_{-\pi}^\pi u(re^{i\theta})d\theta-u(0)=\int_0^r\frac{n(t)}{t}dt,$$
where $n(t)=\mu(\{ z:|z|\leq t\})$, and $\mu$
is the Riesz measure  of $u$, that is the measure with
the density
\begin{equation}\label{3}
\frac{1}{2\pi}\Delta u=\frac{1}{\pi}\| f'\|^2.
\end{equation}
Now positivity of $u$ and (\ref{0}) imply that
all $f_j$ are of order at most $2$, normal type.

In particular,
$$f_j=e^{P_j},\quad 1\leq j\leq n,$$
where
$P_j$ are polynomials of degree at most two.

First we state a lemma which is
the core of our arguments. It is a refined version
of Lemma 1 in \cite{CH}. We denote by $B(a,r)$
the open disc of radius $r$ centered at the point $a$.
\vspace{.1in}

\noindent
{\bf Lemma 1.} {\em Let $u$ be a non-negative
harmonic function in the closure of the
disc $B(a,R)$, and assume that
$u(z_1)=0$ for some point $z_1\in\partial B(a,R).$
Then 
$$|\Del u(z_1)|\geq\frac{u(a)}{2R}.$$
}

{\em Proof.} The function
$$b(r)=\min_{|z-a|=r} u(z)$$
is decreasing and $b(R)=0$.
Harnack's inequality gives
$$b(t)\geq \frac{R-t}{R+t}u(a),\quad 0\leq t\leq R.$$
As $$b(t)=|b(R)-b(t)|\leq(R-t)\max_{[t,R]}|b'|,$$
we conclude that for every $t\in(0,R)$ there exists
$r\in[t,R]$ such that
$$|b'(r)|\geq \frac{1}{R-t}\frac{R-t}{R+t}u(a)=
\frac{u(a)}{R+t}.$$
According to Hadamard's three circle theorem,
$rb'(r)$ is a negative decreasing function, so
$$|Rb'(R)|\geq|rb'(r)|
\geq r\frac{u(a)}{R+t}\geq t\frac{u(a)}{R+t},$$
and the last expression tends to $u(a)/2$ as $t\to R$.
So we have $|b'(R)|\geq u(a)/(2R).$
On the other hand, 
$\displaystyle|\Del u(z_1)|\geq \left|\frac{du}{dn}(z_1)\right|\geq 
|b'(R)|,$
where $d/dn$ is the normal derivative. 
This completes the proof. 
\vspace{.2in}

\noindent
{\bf Proof of the theorem}
\vspace{.2in}

We may assume without loss of generality
that $f_0$ has at least one zero.
Indeed, we can compose $f$ with an automorphism
of $\P^n$, for example replace $f_0$ by $f_0+cf_1,\; c\in\C$
and leave all other $f_j$ unchanged.
This transformation changes neither the $n$
omitted hyperplanes nor the rate of growth of $T(r,f)$
and multiplies the spherical derivative by
a bounded factor. 

Put $u_j=\log|f_j|$, and 
$$u^*=\max_{1\leq j\leq n}u_j.$$
Here and in what follows $\max$ denotes the
pointwise maximum of subharmonic functions. 
We are going to prove first that
\begin{equation}
\label{main}
u_0(z)\leq u^*(z)+4(n+1)|z|\sup_{\C}\| f'\|.
\end{equation}
for $|z|$ sufficiently large.

Let $a$ be a point such that $u_0(a)>u^*(a)$.
Consider the maximal disc  $B(a,R)$ centered at $a$
where the inequality $u_0(z)>u^*(z)$ still holds.
If $z_0$ is a zero of $f_0$ then
$u_0(z_0)=-\infty$ and we have
\begin{equation}\label{nn}
R\leq |a|+|z_0|\leq 2|a|, 
\end{equation}
for $|a|>|z_0|$.
There is a point $z_1\in\partial B(a,R)$ and an
integer $k\in\{1,\ldots,n\}$ such that
\begin{equation}\label{order}
u_0(z_1)=u^*(z_1)=u_k(z_1)\geq u_j(z_1),
\end{equation}
for all $j\in\{1,\ldots,n\}$.
Applying Lemma 1 to the positive harmonic function
$u_0-u_k$ in $B(a,R)$ we obtain
$$\left|\nabla(u_0-u_k)(z_1)\right|\geq
\frac{u_0(a)-u_k(a)}{2R},$$
or
\begin{equation}\label{1n}
u_0(a)\leq u_k(a)+2R\left|\nabla u_0(z_1)-\nabla u_k
(z_1)\right|.
\end{equation}
On the other hand, $|f_0(z_1)|=|f_k(z_1)|\geq|f_j(z_1)|$
for all $j\in\{1,\ldots,n\}$, so
\begin{equation}\label{2n}
\| f'(z_1)\|\geq
\frac{|f_0^\prime(z_1)f_k(z_1)
-f_0(z_1)f_k^\prime(z_1)|}{|f_0(z_1)|^2+
\ldots+|f_n(z_1)|^2}%\nonumber
\geq(n+1)^{-1}\left|\frac{f_0^\prime(z_1)}{f_0(z_1)}-
\frac{f_k^\prime(z_1)}{f_k(z_1)}\right|.
\end{equation}
Combining (\ref{1n}), (\ref{2n}) and (\ref{nn}),
and taking into account
that $|\nabla\log|f||=|f'/f|$, we obtain (\ref{main}).

If all polynomials $P_j$ are linear then inequality
(\ref{main}) completes the proof. Suppose now that
some $P_j$ is of degree $2$. 

Consider again the subharmonic functions $u_j=\log|f_j|,\;
0\leq j\leq n$. For each $j\in\{0,\ldots,n\}$,
the family
$$\{ r^{-2}u_j(rz): r>0\}$$ 
in uniformly bounded from above on compact subsets
of the plane, and bounded from below at $0$. By 
\cite[Theorem 4.1.9]{Hor} these families are normal
(from every sequence one can choose a subsequence that
converges in $L^1_{\rm loc}$).
Take a sequence $r_k$ such that 
\begin{equation}
\label{2}
\lim_{k\to\infty}
\frac{1}{r_k^2}\int_{-\pi}^\pi u(r_ke^{i\theta})d\theta>0,
\end{equation}
where $u$ is defined in (\ref{u}).
Such sequence exists because we assume that at least
one of the $P_j$
is of degree two. 

Then we choose a subsequence (still denoted by $r_k$)
such that 
$$r_k^{-2}u_j(r_kz)\to v_j,\quad 0\leq j\leq n,$$
and $r_k^{-2}u(r_kz)\to v$, where $v_j,v$ are some
subharmonic functions in $\C$.
Then
$$v=\max\{ v_0,\ldots,v_n\}\neq 0$$
is a non-negative subharmonic function.
Let $\nu$ be the Riesz measure of $v$.
 Notice that $\nu\neq 0$
because $v$ is non-negative and $v\neq 0$.
We have weak convergence
$$\nu=\lim_{k\to\infty}\mu_{r_k},$$
where $$\mu_{r_k}(E)=r_k^{-2}\mu(r_kE)$$ for every
Borel set $E$.
Now (\ref{3}) and the condition that $\| f'\|$ is bounded
imply 
\vspace{.1in}

\noindent
{\bf Lemma 2.} {\em $\nu$ is absolutely continuous
with respect to
Lebesque's measure in the plane, with bounded density.}
\vspace{.1in}

{\em Proof.} For every disc $B(a,\delta)$ we have
$$\nu(B(a,\delta))\leq\liminf_{k\to\infty}r_k^{-2}\mu(r_ka,
r_k\delta)\leq\delta^2\sup_{\C}\| f'\|^2.$$

Now we invoke our inequality (\ref{main}). It implies
that 
$$v_0\leq v^*=\max( v_1,\ldots,v_n),$$ so $v=v^*$.
Thus the measure $\nu$ is
supported by finitely many rays. This contradiction
with Lemma 2 shows that all polynomials $P_j$
are in fact linear.
This completes the proof.
\vspace{.2in}

\noindent
{\bf Example}
\vspace{.2in}

Let $\Gamma_0=\{ n+im: n,m\in\Z\}$ be the integer lattice
in the plane, and $\Gamma_1=\Gamma+(1+i)/2$.
For $j\in\{0,1\},$ let $f_j$ be the Weierstrass
canonical products of genus 2
with simple zeros on $\Gamma_j$. 
Then the $f_j$ are entire functions of completely regular
growth in the sense of Levin--Pfluger and their zeros
satisfy the $R$-condition in \cite[Theorem 5, Ch. 2]{Levin}.
This theorem of Levin implies that
\begin{equation}\label{B0}
\log|f_j(re^{i\theta})|=(c+o(1))r^2,
\end{equation}
as $r\to\infty, \; re^{i\theta}\notin C_0$ where $C_0$
is a union of discs of radius $1/4$ centered at the zeros
of $f_j$.
It follows that
\begin{equation}
\label{B}
|f_0(z)|^2+|f_1(z)|^2\to\infty,\quad z\to\infty.
\end{equation}
Cauchy's estimate for the derivative and (\ref{B0}) give
$$\log|f_j^\prime(z)|\leq (c+o(1))|z|^2,\quad z\to\infty.$$
So for the curve $f=(f_0,f_1,1,\ldots,1)$ we obtain
\begin{eqnarray*}
\| f'\|^2&=&\frac{\sum_{i\neq j}|f_i^\prime f_j-
f_if_j^\prime|^2}{\| f\|^4}\leq
\frac{\left(|f_0^\prime f_1-f_0f_1^\prime|^2+
n(|f_0^\prime|^2+|f_1^\prime|^2)\right)}{(|f_0|^2+|f_1|^2)^2}\\
&=&\frac{|g'|^2}{(1+|g|^2)^2}
+o(1).
\end{eqnarray*}
The spherical derivative of $g$ is bounded because $g$
is an elliptic function. Thus
$f$ is a Brody curve that omits $n-1$ hyperplanes
in general position. Evidently $T(r,f)\sim c_1r^2$.

{\em Purdue University, West Lafayette IN 47907 USA

eremenko@math.purdue.edu}
\end{document}